\documentclass[12pt]{article}
\usepackage{amssymb}
\def\appendix{\par}  % to see appendix
\def\om{\omega}
\def\al{\alpha}
\def\ga{\gamma}
\def\uu{{\mathcal U}}
\def\vv{{\mathcal V}}
\def\ff{{\mathcal F}}
\def\gg{{\mathcal G}}
\def\rr{{\mathbb R}}
\def\infsets{[\om]^\om}
\def\sq{\subseteq}
\def\sm{\setminus}
\def\bb{{\mathfrak b}}
\def\dd{{\mathfrak d}}
\def\pp{{\mathfrak p}}
\def\cc{{\mathfrak c}}
\def\rmand{{\mbox{ and }}}
\def\rmiff{{\mbox{ iff }}}
\def\sigc{{\bf \Sigma}^0_2}
\def\siga{{\bf \Sigma}^1_1}
\def\poset{{\mathbb P}}
\def\forces{{| \kern -2pt \vdash}}

\def\sqbdd{{\sqsubseteq^{{\rm bounded}}}}
\newtheorem{theorem}{Theorem}
\newtheorem{lemma}[theorem]{Lemma}
\newtheorem{question}[theorem]{Question}
\newtheorem{prop}[theorem]{Proposition}
\def\proof{\par\noindent Proof\par\noindent}
\def\qed{\par\noindent QED\par}

\begin{document}

\begin{center}
The cardinal characteristic for relative $\ga$-sets
\end{center}

\begin{flushright}
Arnold W. Miller \footnote{
Thanks for partial support from the BEST conference 2003,
support from Boise State University, and to 
Tomek Bartoszynski, Justin Moore, and Marion Scheepers, for their
hospitality during the time the main result in this paper was 
obtained.
}.
\end{flushright}

\begin{quote}
Abstract:  For $X$ a separable metric space define 
$\pp(X)$ to be the smallest cardinality of
a subset $Z$ of $X$ which is not a relative $\ga$-set in $X$, i.e.,
there exists an $\om$-cover of $X$ with no $\ga$-subcover of
$Z$.  We give a characterization of $\pp(2^\om)$ and
$\pp(\om^\om)$ in terms of definable free filters on $\om$
which is related to the psuedointersection number $\pp$.  
We show that for every uncountable standard analytic space $X$ that either
$\pp(X)=\pp(2^\om)$ or $\pp(X)=\pp(\om^\om)$.  We show that 
both of following statements are each relatively consistent with ZFC:
 (a) $\pp=\pp(\om^\om) < \pp(2^\om)$ and
(b) $\pp < \pp(\om^\om) =\pp(2^\om)$
\end{quote}

First we define $\ga$-set.  An open cover $\uu$ of a  topological space $X$
is an $\om$-cover iff for every finite $F\sq X$ there exists $U\in \uu$ with
$F\sq U$. The space $X$ is a $\ga$-set iff for every $\om$-cover of $X$ there
exists a sequence $(U_n\in\uu:n<\om)$ such that for every $x\in X$ for all
but finitely many $n$ we have $x\in U_n$, equivalently
$$X=\bigcup_{m<\om}\bigcap_{n>m}\;U_n \;\;\;\mbox{ or }\;\;\;
\forall x\in X\;\forall^{\infty} n\in\om\;\; x\in U_n.$$
We refer to the sequence $(U_n:n<\om)$ as a
$\ga$-cover of $X$, although technically we are supposed to assume
that the $U_n$ are distinct.
In this paper all our spaces are
separable metric spaces, so we may assume that all $\om$-covers are
countable.  This is because we can replace an arbitrary  $\om$-cover with a
refinement consisting of finite unions of basic open sets.  

The $\ga$-sets were first considered by Gerlits and Nagy \cite{gn}.
One of the things that they showed was the following.
The psuedointersection number $\pp$ is defined as follows:
 $$\pp=\min\{|\ff|\;:\;\ff\sq\infsets\mbox{ has the FIP and }\neg\exists
 X\in\infsets\; \forall Y\in\ff\; X\sq^* Y\}$$
where FIP stands for the finite intersection property, i.e., every
finite subset of $\ff$ has infinite intersection, and $\sq^*$ denotes
inclusion mod finite.  The set $X$ in this definition is called the 
pseudointersection of the family $\ff$.

Gerlits and Nagy \cite{gn} showed that every $\ga$-set 
has strong measure zero (in fact, the Rothberger property
$C^{\prime\prime}$)
and that Martin's Axiom implies every set of reals of size smaller than the
continuum is a $\ga$-set.  Their arguments show that
$$\pp={\rm non}(\ga\mbox{-set })=^{def}\min\{|X|\;:\; X \mbox{ is not a
 $\ga$-set }\}$$
where we only consider separable metric spaces $X$.

The property of being a $\ga$-set is not hereditary.
In fact, a $\ga$-set $X$ of size continuum is constructed
in Galvin and Miller \cite{gm} using MA, which has the property that
there exists a countable $F\sq X$ such that
$X\sm F$ is not a $\ga$-set.  However, any closed subspace of a $\ga$-set  
is a $\ga$-set.

Babinkostova, Guido and Kocinac \cite{babink} have defined the notion of a 
relative $\ga$-set. 
This is also studied in Babinkostova, Kocinac, and Scheepers \cite{babink2}.
For $X\sq Y$ define $X$ is a $\ga$-set relative to $Y$ iff
for every open $\om$-cover $\uu$ of $Y$ there exists
a sequence $(U_n\in\uu:n<\om)$ such that 
 $$X \sq \bigcup_{m<\om}\bigcap_{n>m}\;U_n.$$
Note that if $Z\sq X \sq Y$ and $X$ is a relative $\ga$-set in
$Y$, then $Z$ is also.  

Define the following cardinal number:
$$\pp(Y)=\min\{|X|\;:\; X\sq Y \mbox{ is not a $\ga$-set relative to } Y\}.$$
Perhaps it should be written non($\ga$ relative to $Y$).

In Just, Scheepers, Szeptycki, and Miller \cite{jmss} many cardinal
characteristics for covering properties are shown to be equal to well-known
cardinals.    Scheepers has noted that the cardinal
numbers of the relativized version of the Rothberger property 
$C^{\prime\prime}$ work out to be
either cov(meager) (the cardinality of the smallest cover of the real line with
meager sets) or non(SMZ) (the cardinality of the smallest non strong measure
zero set of reals).

Scheepers has raised the question of what we can say about the
relativized versions for the $\ga$-property.  We begin with the easy

\begin{prop}
$\pp\leq \pp(\om^\om)\leq \pp(2^\om)\leq \cc$
\end{prop}
\proof
If $X$ is a $\ga$-set, then it is a $\ga$-set relative to 
any superspace. Let $|X|=\pp(\om^\om)$ be a subset of $\om^\om$
which is not a relative $\ga$-set.  Then $X$ is not a $\ga$-set relative
to itself, and hence $\pp\leq |X|=\pp(\om^\om)$.

For the second inequality, suppose $X\sq 2^\om$ is not a 
$\ga$-set relative to $2^\om$ with $|X|=\pp(2^\om)$. 
Let $\uu$ be an $\om$-cover of $2^\om$ witnessing that $X$ is
not a relative $\ga$-set in $2^\om$.   
Then
$$\{U\cup (\om^\om\sm 2^\om)\;:\; U\in\uu\}$$
is an $\om$-cover of $\om^\om$ witnessing that $X$ 
is not a $\ga$-set relative to $\om^\om$, and so
$\pp(\om^\om)\leq |X|=\pp(2^\om)$.
\qed

We give another characterization of $\pp(\om^\om)$ and
$\pp(2^\om)$. A filter is free iff it contains the cofinite sets.
For $\ff\sq P(\om)$ a free filter on $\om$, define
$$\pp_\ff=\min\{|X|\;\: X\sq\ff\;\;
\neg\exists a\in\infsets\; \forall b\in X\; a\sq^* b\}.$$
Note that $\pp$ is the minimum of $\pp_\ff$ for
$\ff\sq P(\om)$ a free filter, since every family with the FIP generates
a filter.
We have the following characterizations:

\begin{theorem}\label{characterization}
\par (a) $\pp(\om^\om)$ is the minimum
of $\pp_\ff$ such that $\ff\subseteq P(\om)$ is a ${\bf \Sigma}^1_1$
free filter.
\par (b) $\pp(2^\om)$ is the minimum of
$\pp_\ff$ such that $\ff\sq P(\om)$ is a $\sigc$ free filter.
\end{theorem}
\proof

Suppose $X\sq \om^\om$ with $|X|=\pp(\om^\om)$ and 
$\uu$ is an open $\om$-cover of $\om^\om$ witnessing that $X$ is
not a relative $\ga$-set.  
Without loss of generality we may assume that $\uu$ is a countable
family of clopen sets, say $\uu=\{U_n:n\in\om\}$.  Let $f:\om^\om\to P(\om)$
be the Marczewski \cite{marc} characteristic function of sequence
  $$f(x)=\{n:x\in U_n\}.$$  
This is a continuous mapping so
its image is ${\bf \Sigma}^1_1$.  Since $\uu$ was an $\om$-cover the
image has
the FIP and note that the filter $\ff$ generated by a
${\bf \Sigma}^1_1$
family $\gg$ with the FIP is ${\bf \Sigma}^1_1$, i.e.,
$$X\in\ff \rmiff \exists F\in [\gg]^{<\om}\;\; \cap F\sq X.$$
Now assume $|X| < \pp(\ff)$.  Then
there exists $a\in\om^\om$ such that for each $b\in X$ we have that
$a\sq^* f(b)$.  It follows that $\langle U_n\;:\;n\in a\rangle$ is a $\ga$-cover of
$X$ which is a contradiction.  Hence $\pp(\om^\om)=|X|\geq \pp_\ff$ and so

$$\pp(\om^\om)\geq \min\{\pp_\ff:\ff \mbox{ is a $\siga$ free filter }\}.$$

To see that the other inequality suppose $\ff$ is a ${\bf \Sigma}^1_1$
filter and $X\sq \ff$ witnesses
$|X|=\pp_\ff$.  Let $f:\om^\om\to\ff$ be a continuous onto map.  For each
$n\in\om$ define $U_n=f^{-1}(\{x\in\ff\;:\; n\in x\})$.  Define
$\uu=\{U_n\;:\;n\in\om\}$.  Then $\uu$ is an $\om$-cover of $\om^\om$.
Choose $Y\sq\om^\om$ with $f(Y)=X$ and $|Y|=|X|$. 
If $Y$ is relative $\ga$ in $\om^\om$,
then there exists $a\in\infsets$
such that $\langle U_n\;:\;n\in a\rangle$ is a
$\ga$-cover of $Y$.  For each $b\in X$ we have $c\in Y$ with $f(c)=b$.
For each $n$ if $c\in U_n$, then
$f(c)\in f(U_n)$ and so $n\in b$.  It follows that $a\sq^*c$
for all $c\in X$. Since we are assuming that there is no such $a$, we
must have that $Y$ is not a $\ga$-set relative to $\om^\om$ and therefore
$$\pp(\om^\om)\leq |Y|=|X|=\min\{\pp_\ff:\ff \mbox{ is a
 $\siga$ free filter }\}.$$

The proof for $\pp(2^\om)$ is similar.
Suppose $X\sq 2^\om$ and $\uu$ is an countable clopen $\om$-cover of $2^\om$.
Let $f:2^\om\to P(\om)$
be $f(x)=\{n:x\in U_n\}$. Then $f$ is continuous and so
its range is a compact subset $C\sq P(\om)$ which has the FIP.  For
each $n<\om$ let $h:C^n\to P(\om)$ be defined by 
$$h(X_1,\ldots,X_n)=X_1\cap X_2\cap\cdots\cap X_n.$$
Then since $h$ is continuous its range $C_n$ is compact.  Also
since the projection of the compact set
  $$\{(x,y)\;:\;x\in C_n \rmand x\sq y\}$$ 
onto the second coordinate is 
compact, we see that the filter $\ff$ generated by $C$ is $\sigc$ in
$P(\om)$.  Hence, if $X\sq 2^\om$ has $|X|<\pp_\ff$, then
there exists $a\in \infsets$ with $a\sq^* f(x)$ for each $x\in X$
and there for $x\in U_n$ for all but finitely many $n\in a$.

Conversely, suppose that $\ff$ is a $\sigc$ free filter in $P(\om)$.  
Then there exists a compact $C\sq\ff$ such that for every $x\in\ff$
there exists a $y\in C$ with $x=^*y$.  To see this, suppose that
$\ff=\cup_{n<\om}C_n$.  For each $n<\om$ let
$$C_n^*=\{x\sq\om : n\sq x\rmand
\exists y\in C_n\;\; \forall i\geq n (i\in y\rmiff i\in x)\}$$
then $C=\cup_{n<\om} C_n^*$ does the trick.   Now suppose 
$X\sq \ff$ with $|X|<\pp(2^\om)$.   Let $f:2^\om\to C$ be continuous
and onto and choose $Y\sq 2^\om$ with $|Y|=|X|$ such that
for each $x\in X$ there exists $y\in Y$ with $f(y)=^*x$.
Let $U_n=f^{-1}(\{x\in C\;:\; n\in x\})$.  Then $\uu=\{U_n:n<\om\}$
is an $\om$-cover of $2^\om$ and so there exists 
$a\in \infsets$ such that for every $y\in Y$ we have that
$y\in U_n$ for all but finitely many $n\in a$ and hence
for each $x\in X$ there is $y\in Y$ with $a\sq^* f(y)=^*x$.

\qed

For another paper studying the connection between $\ga$-sets and
free filters, see LaFlamme and Scheepers \cite{ls}.

\begin{lemma}\label{lem1}
\par (a) Suppose that $X$ is homeomorphic to a closed subspace of
$Y$, then $\pp(Y)\leq \pp(X)$.
\par (b) Suppose that $f:X\to Y$ is continuous and onto, then 
$\pp(X)\leq \pp(Y)$.
\end{lemma}
\proof

\medskip\noindent (a) Suppose $Z\sq X$ with $|Z|=\pp(X)$ is not relatively 
$\ga$ in $X$ and this is witnessed by a family of open sets of $Y$ which
is an $\om$-cover of $X$. 
Then
$$\{U\cup(Y\sm X)\;:\; U\in\uu\}$$
is an $\om$-cover of $Y$ which shows that $Z$ is not relatively $\ga$ in
$Y$.  Hence $\pp(Y)\leq |Z|=\pp(X)$.

\medskip\noindent (b) Suppose $Z\sq Y$ with $|Z|=\pp(Y)$ is not relatively 
$\ga$ in $Y$ and this is witnessed by an $\om$-cover
$\uu$.  Choose $W\sq X$ with $|W|=|Z|$ and $f(W)=Z$.  
Let $\vv=\{f^{-1}(U)\;:\; U\in\uu\}$.  Since $f$ is onto, $\vv$ is
an $\om$-cover of $X$.  We claim that there is no sequence
$(f^{-1}(U_n): n<\om,U_n\in\uu)$ such that for every $x\in W$ for
all but finitely many $n$ we have $x\in f^{-1}(U_n)$.  This is
because $x\in f^{-1}(U_n)$ implies $f(x)\in U_n$, but then
$f(W)=Z$ would have the property that every $y\in Z$ is in all
but finitely many $U_n$.  It follows that $\pp(X)\leq |W|=|Z|=\pp(Y)$.
\qed

\begin{theorem}
Suppose $X$ is an uncountable $\siga$ set in a Polish space,
i.e., a nontrivial standard analytic space, then
\par (a) if $X$ is not $\sigma$-compact, then $\pp(X)=\pp(\om^\om)$ and
\par (b) if $X$ is $\sigma$-compact, then  $\pp(X)=\pp(2^\om)$.
\end{theorem}
\proof

If $X$ is , it contains a homeomorphic copy of $2^\om$ 
and is the continuous
image of $\om^\om$.  It follows from Lemma \ref{lem1} that
$$\pp(\om^\om)\leq \pp(X)\leq \pp(2^\om).$$

\medskip\noindent (a) If $X$ is not $\sigma$-compact, then Hurewicz
\cite{hurewicz} (see Kechris \cite{kechris} 21.18) proved that there exists a
closed subspace of $X$ which is homeomorphic to $\om^\om$.  Hence by Lemma
\ref{lem1} we have $\pp(X)\leq\pp(\om^\om)$.  

\medskip\noindent (b) First suppose that $X=\om \times 2^\om$.
Let $Y\sq \om \times 2^\om$ be 
non relatively $\ga$ with $|Y|=\pp(\om \times 2^\om)$.
Since it is zero dimensional
we can assume that this is witnessed by a countable $\om$-cover 
of clopen sets $\uu=\{C_n:n<\om\}$.  
As in the proof of Theorem \ref{characterization}
we take $f:\om \times 2^\om\to P(\om)$ 
defined by 
  $$f(x)=\{n<\om: x\in C_n\}.$$  
The function $f$ 
is continuous since the $C_n$ are clopen and its image $f(\om\times 2^\om)$
is a 
$\sigma$-compact family of sets with the finite intersection property.
Then $f(\om\times 2^\om)$ generates a $\sigma$-compact filter $\ff$ and 
$f(Y)$ is a subset of $\ff$ without a pseudo-intersection.  Hence
$\pp_\ff\leq |f(Y)|\leq |Y|=\pp(\om \times 2^\om)$ and so we have 
$\pp(2^\om)\leq \pp(\om \times 2^\om)$.

Now suppose that $X$ is any $\sigma$-compact metric space.
Note that there is a continuous onto mapping $f:\om \times 2^\om \to X$
and so by Lemma \ref{lem1} we have that 
$$\pp(X)\geq \pp(\om \times 2^\om)=\pp(2^\om).$$
 
\qed

The main result of this paper is the following theorem:

\begin{theorem}\label{main}
Both of following statements are each relatively consistent with ZFC:
\par (a) $\pp=\pp(\om^\om) < \pp(2^\om)$ and
\par (b) $\pp < \pp(\om^\om) =\pp(2^\om)$
\end{theorem}
\proof

\medskip Part(a). The forcing notion we use is the obvious one.  The only 
difficulty of the proof is that the forcing does only what we want it to
and not more.

Given an $\om$-cover $\uu$ of $2^\om$ 
define  the poset $\poset(\uu)$ as follows:

\begin{enumerate}
\item $p\in \poset(\uu)$  iff $p=(F,(U_n\in\uu:n<N))$ where $N<\om$ and
$F\in [2^\om]^{<\om}$.  

\item $p\leq q$ iff $F^p\supseteq F^q$, $N^p\geq N^q$,
$U_n^p=U_n^q$ for each $n<N^q$, and 
$x\in U_n^p$  for each $x\in F^q$ and $n$ with $N^q\leq n < N^p$.
\end{enumerate}

This poset is the obvious one for generically creating a
$\ga$-subcover of $\uu$ for the ground model elements of $2^\om$.

\begin{lemma}\label{routine}
The partial order $\poset(\uu)$ is $\sigma$-centered. Furthermore,
suppose $G$ is $\poset(\uu)$-generic
over $V$.  Define $(U_n:n<\om)$ by $U_n=U_n^p$ for any $p\in G$ with
$N^p>n$.  Then $\forall x\in V\cap 2^\om\;\;\forall^\infty n\;\;
x\in U_n$.
\end{lemma}
\proof
$\sigma$-centered is clear, since if $(N^p_n:n<N^p)=(N^q_n:n<N^q)$
then the condition $(F^p\cup F^q,(N^p_n:n<N^p))$ extends both $p$ and
$q$. The fact that $U_n$ is defined for every $n<\om$ follows from
$\uu$ being an $\om$-cover and a density argument, i.e.,
given any $p$ with $N_p\leq n$ extend it by adding $U_k$ which
cover $F_p$.  That 
$(U_n:n<\om)$ is a $\ga$-cover of $2^\om\cap V$
since given any $x\in 2^\om$ in the ground model, the set
$$D=\{p\in\poset(\uu)\;:x\in F_p\}$$
is dense and if $x\in F_p$ for some $p\in G$ then $x\in U_n$ 
for every $n\geq N_p$.
\qed

The model for $\pp=\pp(\om^\om) < \pp(2^\om)$ is obtained by starting
with a model of GCH and doing a finite support iteration of
$\poset(\uu_\al)$ for $\al<\om_2$ where at each stage in the iteration
 $$V[G_\al]\models \uu_\al\mbox{ is an $\om$-cover of } 2^\om$$ 
and where we have dove-tailed so as to ensure that for any $\uu$
such that
 $$V[G_{\om_1}]\models \uu\mbox{ is a countable $\om$-cover of } 2^\om$$
then for some $\al<\om_2$ we have that $\uu=\uu_\al$.  This dovetailing
can be done since there are only continuum many 
countable $\om$-covers of $2^\om$ and the intermediate models satisfy
the continuum hypothesis.
 In the model $V[G_{\om_2}]$ we have that $\pp(2^\om)=\omega_2$, so
we need only show that $\pp(\om^\om)=\om_1$.  As usual, define
Rothberger's unbounded number:
$$\bb=\min\{|X|\;:\; X\sq\om^\om\;\; \forall g\in\om^\om\;\; \exists f\in X\;\;
\exists^\infty n\;\; f(n)>g(n)\}.$$

\begin{lemma}\label{inequal}
$\pp(\om^\om)\leq \bb$
\end{lemma}
\proof
Suppose $X\sq \om^\om$ and $|X|<\pp(\om^\om)$.  We need
to show that $X$ is eventually dominated.  Without loss of
generality we may assume that the elements of $X$ are increasing
and $X$ is infinite.
For each $n<\om$ let
$$\uu_n=\{U_m^n:m<\om\}\mbox{ where } U_m^n=\{f\in\om^\om: f(n)<m\}.$$
Each $\uu_n$ is an $\om$-cover of $\om^\om$.
There is a standard trick due to Gerlits and Nagy \cite{gn} for replacing
a sequence of $\om$-covers by a single $\om$-cover.  Let
  $$\{x_n:n<\om\}\sq X$$
be distinct and let
  $$\uu=\{U\sm\{x_n\}:n<\om,\; U\in\uu_n\}.$$ 
Then $\uu$ is an $\om$-cover of $\om^\om$, since given any
finite set $F$ then $x_n\notin F$ for large enough $n$ and
so $F\sq U\sm \{x_n\}$ for some $U\in\uu_n$. 

Since $X$ is a relative $\ga$-set in $\om^\om$ there exists a
sequence from $\uu$ which is a $\ga$-cover of $X$.  
Now since we threw out $x_n$ from each element $\uu_n$ at most finitely
many of the elements of this sequence can come from the same $\uu_n$.  
By taking an infinite subsequence we may assume that
$(U_{g(n)}^n:n\in A)$ is a $\ga$-cover of $X$ for some infinite
$A\sq\om$.
It follows that 
for every $f\in X$ that 
$$\forall^\infty\; n\in A\;\;f(n)<g(n).$$ 
Since the $f\in X$ are increasing if
we extend $g$ to all of $\om$ by letting $g(m)=g(n)$ where
$n\in A$ is minimal so that $n\geq m$, then $g$ eventually dominates
every $f\in X$ on all of $\om$.  

It follows that $|X|<\bb$.  Since $X$ was arbitrary
we get that $\pp(\om^\om)\leq \bb$.
\qed

Our goal is to show that $\bb=\om_1$ holds in this model.
For the next two lemmas we assume $\uu$ is an $\om$-cover of $2^\om$
and the forcing is $\poset(\uu)$.

\begin{lemma} \label{finite} 
Suppose we are given $(U_n\in\uu:n<N)$, $k<\om$, and a term $\tau$ such
that $\forces \tau\in\om$.  Then there exists $K<\om$ such that
for every $p\in\poset(\uu)$ with $|F^p|\leq k$ and
$(U_n\in\uu:n<N)=(U_n^p\in\uu:n<N^p)$ there exists $q\leq p$ such that
$q\forces \tau<K$.
\end{lemma}
\proof
Call $q\in\poset(\uu)$ good iff
\begin{enumerate}
\item $N^q\geq N$
\item $U_n=U_n^q$ for all $n<N$, and
\item $q$ decides $\tau$, i.e. for some $m$, $\;\;q\forces \tau=m$.
\end{enumerate}
For good $q$ define:
$$U_q=\{(x_1,\ldots,x_k)\in (2^\om)^k\;:\; \forall i\;\;\;
 (N\leq i < N_q)\to 
\{x_1,\ldots,x_k\}\sq U_i^q\}.$$
Note that each $U_q$ is an open subset of $(2^\om)^k$.  
Also the family $U_q$ for $q$ good cover 
$(2^\om)^k$. This is because given any $(x_1,\ldots,x_k)$ there
exist a condition $q\leq (\{x_1,\ldots,x_k\},(U_n:n<N))$ which
decides $\tau$ and therefor is good.  By compactness there exist
finitely many good $q$, say $\Gamma$, such that
$\{U_q:q\in\Gamma\}$ covers $(2^\om)^k$.

Since each good $q$
decides $\tau$ and $\Gamma$ is finite, 
we can find $K$ so that for each $q\in \Gamma$
$$q\forces \tau<K.$$ 
Note
that for any $p$ as in the Lemma, if $F^p\sq U_q$ then $q$ and $p$ are
compatible since $(F^p\cup F^q,(N^q_n:n<N^q))$ extends both of them.
\qed

It is not hard to see from this lemma that our forcing does not add a
dominating sequence.   In order to prove the full result we need to show this
for the iteration. To do this we prove the following
stronger, but more technical, property 
(see Bartoszynski and Judah \cite{barto} definition
6.4.4).

\begin{lemma}
The poset $\poset(\uu)$ is really $\sqbdd$-good, i.e., for every name
$\tau$ for an element of $\om^\om$ there exists $g\in\om^\om$ such 
that for any $x\in\om^\om$ if
there exists $p\in\poset(\uu)$ such that 
$p\forces$``$\;\forall^\infty n\; x(n)<\tau(n)$'',
then $\forall^\infty n\;\; x(n)<g(n)$. 
\end{lemma}
\proof
Let $k_n,(U^n_m\;:\;m<N_n)$ for $n<\om$ 
list with infinitely many repetitions all
pairs of $k<\om$ and finite sequences from $\uu$.  Using Lemma \ref{finite}
repeatedly we can construct $g\in\om^\om$ such that for every 
$l<\om$:

\noindent for any $n < l$ and $p\in\poset(\uu)$ with 
$$|F^p|\leq k_n \rmand (U_m^n:m<N_n)=(U_m^p:m<N^p)$$
there exists $q\leq p$ such that $q\forces\tau(l)<g(l)$.  

Now suppose $p\forces \forall^\infty n\;x(n)<\tau(n)$.  By extending $p$ 
(if necessary) we may assume there exists $n_0$ such that
$p\forces \forall n> n_0\; x(n)<\tau(n)$.  By making $n_0$ larger
(if necessary) we may assume that  
$$|F^p|=k_{n_0}\rmand (U_m^{n_0}:m<N_{n_0})=(U_m^p:m<N^p).$$

\bigskip
\noindent Claim $\forall n>n_0\;\; x(n)<g(n)$.

proof:
Suppose not and $x(l)\geq g(l)$ for some $l>n_0$.  
By our construction of $g$ we have that there exists $q\leq p$ such that
$q\forces \tau(l)<g(l)$. But this means that $q\forces \tau(l)<x(l)$ which
contradicts the fact that $p\forces \forall n> n_0 \tau(n)>x(n)$.
This proves the Claim and the Lemma.
\qed

It follows (see Bartoszynski and Judah \cite{barto} Theorem 6.5.4)
that the finite support iteration using $\poset(\uu_\al)$ at
stage $\al$
does not add a dominating real and so over a ground model which
satisfies CH we have that
$V[G_{\om_2}]$ satisfies that $\bb=\om_1$ and hence $\pp(\om^\om)=\om_1$
by Lemma \ref{inequal}.
This proves Theorem \ref{main} part (a), the consistency of
$\pp=\pp(\om^\om) < \pp(2^\om)$.

\bigskip
Part (b) (the consistency of $\pp < \pp(\om^\om) =\pp(2^\om)$) is simpler.
It is well known that $\pp>\om_1$ implies that $2^{\om_1}=2^\om$.
For example, see Rothberger \cite{roth}.  Now starting with
a ground model $V$ which satisfies $2^\om=\om_2$ and $2^{\om_1}=\om_3$,
do a finite support iteration using
$\poset(\uu_\al)$ at stage $\al<\om_2$ where
$\uu_\al$ is an $\om$-cover of $V[G_\al]\cap \om^\om$.
Dovetail so that $\uu_\al$ for $\al<\om_2$ lists all countable $\om$-covers of
$\om^\om$ in the final model $V[G_{\om_2}]$.
This can be done since in
all these models the continuum is $\om_2$.  
The analogue of Lemma \ref{routine} holds for $\om^\om$ in
place of $2^\om$ so 
in the final model we
have that $\pp(\om^\om)=\om_2$. 
Also we get $\pp=\om_1$ since $2^{\om_1}=\om_3>\om_2=2^\om$.
This finishes the proof of Theorem \ref{main}.
\qed

One obvious question is 

\begin{question}
Is it consistent to have $\pp<\pp(\om^\om)<\pp(2^{\om})$?
\end{question}

\begin{question}(Scheepers)
Are either $\pp(\om^\om)$ or $\pp(2^{\om})$ the same
as some other well-known small cardinal?
See Vaughan \cite{vaughan} for a plethora of such cardinals.
\end{question}

In Laver's model \cite{laver} for the Borel conjecture, we have that
$\bb=\dd=\om_2$ and $\pp(2^\om)=\pp(\om^\om)=\om_1$. 
In Laver's model there is a set of reals of size $\om_1$ which does not have
measure zero, i.e., non(measure)=$\om_1$,  Judah and Shelah \cite{js}, see also
Judah and Bartoszynski \cite{barto} or Pawlikowski \cite{pawl}. But it is
easy to see that $\pp(2^\om)\leq$ non(measure), i.e.,
if $X\sq 2^\om$ and $|X|<\pp(2^\om)$ then $X$ has measure
zero.  Let $\{x_n:n<\om\}\sq X$ be distinct and look at 
   $$\uu=\{C\sq 2^\om\;:\;\exists n\;\; x_n\notin C \mbox{ is clopen  and } 
   \mu(C)<{1\over 2^n}\}.$$
This is an $\om$-cover of $2^\om$ and so there exists a sequence 
$C_n\in\uu$ with $X\sq\cup_n\cap_{m>n}C_m$.  For any $n$ at most finitely
many $C_n$ have measure $>{1\over 2^n}$ 
which shows that $X$ has measure zero.  

It is also true that $\pp(2^\om)\leq$ non(SMZ), i.e., if $|X|<\pp(2^\om)$ then
$X$ has strong measure zero.  The result of  Gerlits and Nagy \cite{gn}, that
$\ga$-sets have the Rothberger property $C^{\prime\prime}$, relativizes to show
that if $X\sq 2^\om$ and $|X|<\pp(2^\om)$, then $X$ has the relative Rothberger
property and this implies that $X$  has strong measure zero.

\begin{question}
Suppose that $Y=\cup_{n<\om}X_n$ is an increasing union where $Y$ is
a separable metric space.  If each $X_n$ is relatively $\ga$ in $Y$, is
$Y$ a $\ga$-set?  If not, suppose each $X_n$ is a $\ga$-set, then
is $Y$ a $\ga$-set?
\end{question}

\begin{flushleft}
Arnold W. Miller \\
miller@math.wisc.edu \\
http://www.math.wisc.edu/$\sim$miller\\
University of Wisconsin-Madison \\
Department of Mathematics, Van Vleck Hall \\
480 Lincoln Drive \\
Madison, Wisconsin 53706-1388 \\
\end{flushleft}

%%%%%%%%%%%%%%%%%%%%%%%%%%%%%%%%%%%%%%%%%%%%%%%%%%%%%%%%%%%%%%%%%%%%%%%

\appendix

\newpage

The appendix is not intended for final publication but for  
the electronic version only.

\begin{center}
Appendix \\
Scheepers Remarks
\end{center}

\noindent Def. $X\sq Y$ is $C^{\prime\prime}$ (Rothberger) in $Y$ iff for every 
sequence $(\uu_n:n<\om)$ of open covers of $Y$ there
exists a cover $(U_n\in\uu_n:n<\om)$ of $X$.

\bigskip\noindent
Prop. \par (a) non($C^{\prime\prime}$ in
$\om^\om$)=non($C^{\prime\prime}$)=cov(meager)=
non(SMZ in $\om^\om$)
\par (b) non($C^{\prime\prime}$ in $2^\om$) = non(SMZ in $2^\om$)=
non(SMZ in $\rr$)

\proof

Here we mean strong measure zero in the usual metric on
the reals and for $\om^\om$ or $2^\om$ the metric $d(x,y)={1\over n+1}$ 
where $n$ is minimal such that $x(n)\not=y(n)$.

\bigskip\noindent
(a) Fremlin and Miller
(1988) prove: 

non($C^{\prime\prime}$)=cov(meager)=non(SMZ in $\om^\om$)

non($C^{\prime\prime})\leq $non($C^{\prime\prime}$ in $\om^\om$)
since if $X$ is not relatively $C^{\prime\prime}$ it is not $C^{\prime\prime}$.

non($C^{\prime\prime}$ in $\om^\om$) $\leq$ non(SMZ in $\om^\om$) since
$C^{\prime\prime}\sq$ SMZ.

\bigskip\noindent
(b)  Suppose $X\sq 2^\om$ fails to
be relatively $C^{\prime\prime}$.
Note that by compactness of $2^\om$ we may assume
there is a sequence $(\uu_n:n<\om)$ of
finite clopen covers of $2^\om$ for which there is no
$(U_n\in\uu_n:n<\om)$ which covers $X$.  Now choose $\epsilon_n>0$ so
that any interval $[s]$ with diameter less than $\epsilon_n$ is
a subset of some $U_n$.  For the converse, suppose $X$ fails to
have SMZ in $2^\om$ which is witnessed by $(\epsilon_n:n<\om)$.
Then the sequence $\uu_n=\{C\sq 2^\om: {\rm diam}(C)<\epsilon_n\}$ witnesses
that it is not $C^{\prime\prime}$.  

\par 
non(SMZ in $2^\om$)=non(SMZ in $[0,1]$)=non(SMZ in $\rr$) is easy to prove.

\qed

These cardinals can be different, for example, in the iterated
modified Silver reals model, see Miller (1981),
cov(meager)$=\om_1$ while non(SMZ in $2^\om)=\om_2$.

\bigskip

\noindent Def. $X$ has the Menger property M iff for every sequence 
$(\uu_n:n<\om)$ of open covers there exists $(\vv_n\in[\uu_n]^{<\om}:n<\om)$
such that $X\sq \cup_n\cup \vv_n$. 

\noindent Def. $X$ has the Hurewicz property H iff for every sequence 
$(\uu_n:n<\om)$ of open covers there exists $(\vv_n\in[\uu_n]^{<\om}:n<\om)$
such that $X\sq \cup_m\cap_{n>m}\cup\vv_n$.

\bigskip
Then it is known, see Miller and Fremlin (1988) that
$$\dd={\rm non}(M)\;\;\;\;\;\bb={\rm non}(H).$$
Relativizing H or M to $2^\om$ doesn't work since $2^\om$ has
property H and M.  For $\om^\om$ it is easy to see:

\bigskip 
Prop.
\par (a) non(M) $\leq$ non(M in $\om^\om$) $\leq\dd\leq $ non(M)
\par (b) non(H) $\leq$ non(H in $\om^\om$) $\leq\bb\leq $ non(H)

\bigskip\bigskip

\noindent Biblio

\bigskip

Miller, Arnold W.; Some properties of measure and category.  Trans. Amer. Math.
Soc.  266  (1981), no. 1, 93--114.

\medskip

Miller, Arnold W.; Fremlin, David H.; On some properties of Hurewicz, Menger,
and Rothberger. Fund. Math. 129 (1988), no. 1, 17--33.

\bigskip\bigskip\bigskip

 TeXed on : \today

\end{document}